\documentclass[runningheads,a4paper]{llncs}

\usepackage{graphicx}
\usepackage[caption=false]{subfig}
\usepackage{amssymb, amsfonts, amsbsy, latexsym, epsfig,
  color,amsmath,pstool}
\usepackage{url}

\urldef{\mailsa}\path|kutyniok,|
\urldef{\mailsb}\path|wlim@uni-osnabrueck.de|
\newcommand{\keywords}[1]{\par\addvspace\baselineskip
\noindent\keywordname\enspace\ignorespaces#1}

\newcommand*{\numbersys}[1]{\ensuremath{\mathbb{#1}}}

\newcommand*{\R}{\numbersys{R}}

\newcommand*{\Z}{\numbersys{Z}}

\usepackage{xspace}

\newcommand\cC{{\mathcal{C}}}

\newcommand\cP{{\mathcal{P}}}

\begin{document}

\mainmatter  

\title{Image Separation using Wavelets and Shearlets}


%
%
\author{Gitta Kutyniok \and Wang-Q Lim}
%

\institute{Institute of Mathematics, University of Osnabr\"uck, \\ 49069 Osnabr\"uck, Germany\\
\mailsa
\mailsb\\}
%
%

\maketitle

\begin{abstract}
In this paper, we present an image separation method for separating images into point- and curvelike
parts by employing a combined dictionary consisting of wavelets and compactly supported shearlets utilizing the fact that
they sparsely represent point and curvilinear singularities, respectively. Our methodology is based
on the very recently introduced mathematical theory of geometric separation, which shows that highly
precise separation of the morphologically distinct features of points and curves can be achieved
by $\ell^1$ minimization. Finally, we present some experimental results showing the effectiveness
of our algorithm, in particular, the ability to accurately separate points from curves even
if the curvature is relatively large due to the excellent localization property of compactly supported
shearlets.
\keywords{Geometric separation, $\ell_1$ minimization, sparse approximation, shearlets, wave\-lets}
\end{abstract}

\section{Introduction}
The task of separating an image into its morphologically different contents has recently drawn a lot of
attention in the research community due to its significance for applications. In neurobiological imaging,
it would, for instance, be desirable to separate 'spines' (pointlike objects) from 'dendrites' (curvelike
objects) in order to analyze them independently aiming to detect characteristics of Alzheimer disease. Also,
in astronomical imaging, astronomers would often like to separate stars from filaments for further analysis,
hence again separating point- from curvelike structures. Successful methodologies for efficiently and
accurately solving this task can in fact be applied to a much broader range of areas in science and
technology including medical imaging, surveillance, and speech processing.

Although the problem of separating morphologically distinct features seems to be intractable -- the problem
is underdetermined, since there is only one known data (the image) and two or more unknowns --
there has been extensive studies on this topic. The book by Meyer \cite{Me02} initiated the
area of image decomposition, in particular, the utilization of variational methods. Some years later,
Starck, Elad, and Donoho suggested a different approach in  \cite{SED05} coined `Morphological Component Analysis', which
proclaims that such a separation task might be possible provided that we have prior information
about the type of features to be extracted and provided that the morphological difference between those
is strong enough. For the separation of point- and curvelike features, it was in fact recently even
theoretically proven in \cite{DK10} that $\ell_1$ minimization solves this task with
arbitrarily high precision exploring a combined dictionary of wavelets and curvelets. Wavelets provide
optimally sparse expansions for pointlike structures, and curvelets provide optimally sparse expansions for
curvelike structures. Thus $\ell_1$ minimization applied to the expansion coefficients of the original image into
this combined dictionary forces the pointlike structures into the wavelet part and the curvelike structures
into the curvelet part, thereby automatically separating the image. An associated algorithmic approach using
wavelets and curvelets has been implemented in MCALab\footnote{MCALab (Version 120) is available from
\url{http://jstarck.free.fr/jstarck/Home.html}.}.

Recently, a novel directional representation system -- so-called shearlets -- has emerged which provides a
unified treatment of continuum models as well as digital models, allowing, for instance, a precise
resolution of wavefront sets, optimally sparse representations of cartoon-like images, and associated fast
decomposition algorithms; see the survey paper \cite{KLL10b}. Shearlet systems are systems generated by one
single generator with parabolic scaling, shearing, and translation operators applied to it, in the same way
wavelet systems are dyadic scalings and translations of a single function, but including a directionality
characteristic owing to the additional shearing operation (and the anisotropic scaling). The shearing operation
in fact provides a more favorable treatment of directions, thereby ensuring a unified treatment of the continuum
and digital realm as opposed to curvelets which are rotation-based in the continuum realm, see \cite{CD04}.

Thus, it is natural to ask whether also a combined dictionary of wavelets and shearlets might be utilizable
for separating point- and curvelike features, the advantage presumably being a faster scheme, a
more precise separation, and a direct applicability of theoretical results achieved for the continuum domain.
And, in fact, the theoretical results from \cite{DK10} based on a model situation were shown to also
hold for a combined dictionary of wavelets and shearlets \cite{DK09}. Moreover, numerical results give
evidence to the superior behavior of shearlet-based decomposition algorithms when compared to curvelet-based algorithms;
see \cite{KLL10b} for a comparison of ShearLab\footnote{ShearLab (Version 1.0) is available from \url{http://www.shearlab.org}.}
with CurveLab\footnote{CurveLab (Version 2.1.2) is available from \url{http://www.curvelet.org}.}.

In this paper, we will present a novel approach to the separation of point- and curvelike features exploiting
a combined dictionary of wavelets and shearlets as well as utilizing block relaxation in a particular way.
Numerical results give evidence that indeed the previously anticipated advantages hold true, i.e., that this
approach is superior to separation algorithms using wavelets and curvelets such as MCALab in various ways, in particular,
our algorithm is faster and provides a more precise separation, in particular, if the curvature of the
curvilinear part is large. In the spirit of reproducible research
\cite{DMSSU09}, our algorithm is included in the freely available ShearLab toolbox.

\medskip

This paper is organized as follows. In Section \ref{sec:shearlets}, we introduce the multiscale system of
shearlets, and Section \ref{sec:geosep} reviews the mathematical theory of geometric separation of point- and
curvelike features. Our novel algorithmic approach is presented in Section \ref{sec:code} with numerical
results discussed in Section \ref{sec:numerics}.

\section{Shearlets}
\label{sec:shearlets}

In most multivariate problems, important features of the considered data are concentrated on lower dimensional
manifolds. For example, in image processing an edge is an 1D curve that follows a path of rapid change
in image intensity. Recently, the novel directional representation system of shearlets \cite{LLKW05,GKL06} has
emerged to provide efficient tools for analyzing the intrinsic geometrical features of a signal using anisotropic
and directional window functions. In this approach, directionality is achieved by applying integer powers of a
shear matrix, and those operations preserve the structure of the integer lattice which is crucial for digital
implementations. In fact, this key idea leads to a unified treatment of the continuum as well as digital realm,
while still providing optimally sparse approximations of anisotropic features. As already mentioned before, shearlet
systems are generated by parabolic scaling, shearing, and translation operators applied to one single generator.
Let us now be more precise and formally introduce shearlet systems in 2D.

We first start with some definitions for later use. For $j \ge 0$ and $k \in \mathbb{Z}$, let
\[
A_{2^j} = \begin{pmatrix} 2^j & 0 \\ 0 & 2^{j/2} \end{pmatrix},
\quad
\tilde{A}_{2^j} = \begin{pmatrix} 2^{j/2} & 0 \\ 0 & 2^{j} \end{pmatrix},
\quad \text{and} \quad
S_k = \begin{pmatrix} 1 & k \\ 0 & 1 \end{pmatrix}.
\]
We can now define so-called cone-adapted discrete shearlet systems, where the term `cone-adapted' originates from
the fact that these systems tile the frequency domain in a cone-like fashion.
For this, let $c$ be a positive constant, which will later control the sampling density. For $\phi, \psi, \tilde{\psi} \in
L^2({\mathbb R}^2)$, the {\em cone-adapted discrete shearlet system} $SH(\phi,\psi,\tilde{\psi};c)$ is then defined by
\[
SH(\phi,\psi,\tilde{\psi};c) = \Phi(\phi;c) \cup \Psi(\psi;c) \cup \tilde{\Psi}(\tilde{\psi};c),
\]
where
\begin{align*}
\Phi(\phi;c) &= \{\phi(\cdot-cm) : m \in {\mathbb Z}^2\},\\
\Psi(\psi;c) &= \{\psi_{j,k,m} = 2^{\frac34 j} \psi(S_kA_{2^j}\, \cdot \, -cm) : j \ge 0, |k| \leq \lceil2^{{j/2}}\rceil, m \in {\mathbb Z}^2\},\\
\tilde{\Psi}(\tilde{\psi};c) &= \{\tilde{\psi}_{j,k,m} = 2^{\frac34 j} \tilde{\psi}(S^T_k\tilde{A}_{2^j}\, \cdot \, -cm) : j \ge 0,
|k| \leq \lceil2^{{j/2}}\rceil, m \in {\mathbb Z}^2 \}.
\end{align*}

In \cite{KLP10}, a comprehensive theory of compactly supported shearlet frames
is provided, i.e., systems with excellent spatial localization.
It should also be mentioned that in \cite{KL10} a large class of compactly supported shearlet frames were shown to provide optimally sparse
approximations of images governed by curvilinear structures, in particular, so-called {\em cartoon-like images}
as defined in \cite{CD04}. This fact will be explored in the sequel.

\section{Mathematical Theory of Geometric Separation}
\label{sec:geosep}

In \cite{SED05}, a novel image separation method -- Morphological Component Analysis (MCA) -- based on
sparse representations of images was introduced. In this approach, it is assumed that each image is the
linear combination of several components that are morphologically distinct -- for instance, points, curves,
and textures. The success of this method relies on the assumption that each of the components is
sparsely represented in a specific representation system. The key idea is then the following: Provided
that such representation systems are identified, the usage of a pursuit algorithm searching for the
sparsest representation of the image with respect to the dictionary combining all those specific
representation systems will lead to the desired separation.

Various experimental results in \cite{SED05} show the effectiveness of this method for image separation however
without any accompanying mathematical justification. Recently, the first author of this paper and Donoho
developed a mathematical framework in \cite{DK10} within which the notion of successful separation can be made
definitionally precise and can be mathematically proven in case of separating point- from curvelike features,
which they coined {\em Geometric Separation}. One key ingredient of their analysis is the consideration of
clustered sparsity properties measured by so-called {\em cluster coherence}.  In this section, we briefly review this theoretical approach to
the Geometric Separation Problem, which will serve as the foundation for our algorithm.

\subsection{Model Situation}
\label{subsec:model}
As a mathematical model for a composition of point- and curvelike structures, we consider the following two
components: As a `point-like' object, we consider the function $\cP$ which is smooth except for point
singularities and is defined by
\[
\cP = \sum_{i=1}^{P}|x-x_i|^{-3/2}.
\]
As a `curve-like' object, we consider the distribution $\cC$ with singularity along a closed curve
$\tau:[0,1] \rightarrow \R^2$ defined by
\[
\cC = \int \delta_{\tau(t)}dt.
\]
Then our model situation is the sum of both, i.e.,
\begin{equation}\label{eq:point_curve}
f = \cP+\cC.
\end{equation}
The Geometric Separation Problem now consists of recovering $\cP$ and $\cC$ from the observed signal $f$.

\subsection{Chosen Dictionary}

As we indicated before, it is now crucial to choose two representations systems each of which sparsely represents
one of the morphologically different components in the Geometric Separation Problem. Our sparse approximation
result, described in the previous section, suggests that curvilinear singularities can be sparsely represented
by shearlets. On the other hand, it is well known that wavelets can provide optimally sparse approximations of
functions which are smooth apart from point singularities. Hence, we choose the overcomplete system within which
we will expand the signal $f$ as a composition of the following two systems:
\begin{itemize}
\item {\it Orthonormal Separable Meyer Wavelets}: Band-limited wavelets which form an orthonormal basis
of isotropic generating elements.
\item {\it Bandlimited Shearlets}: A directional and anisotropic tight frame generated by a band-limited
shearlet generator $\psi$ defined in Section \ref{sec:shearlets}.
\end{itemize}


\subsection{Subband Filtering}
\label{subsec:filtering}

Since the scaling subbands of shearlets and wavelets are similar 
we can define a family of filters $(F_j)_{j}$ which allows to decompose a function $f$
into pieces $f_j$ with different scales $j$ depending on those subbands. The piece $f_j$ associated to subband $j$
arises from filtering $f$ using $F_j$ by
\[
f_j = F_j \ast f,
\]
resulting in a function whose Fourier transform $\hat f_j$ is supported on the scaling subband of scale $j$ of
the wavelet as well as the shearlet frame. The filters are defined in such way, that the original function
can be reconstructed from the sequence $(f_j)_j$ using
\[
f = \sum_{j} F_j \ast f_j, \quad f \in L^2(\R^2).
\]
We can now exploit these tools to attack the Geometric Separation Problem scale-by-scale. For this, we filter the model
problem (\ref{eq:point_curve}) to derive the sequence of filtered images
\[
f_j = \cP_j + \cC_j \quad \text{for all scales} \,\, j.
\]

\subsection{$\ell_1$ Minimization Problem}

Let now $\Phi_1$ and $\Phi_2$ be an orthonormal basis of band-limited wavelets and a tight frame of band-limited shearlets,
respectively. Then, for each scale $j$, we consider the following optimization problem:
\begin{equation}\label{eq:sep}
(\hat W_j,\hat S_j) = \text{argmin}_{W_j,S_j}\|\Phi^T_1 W_j\|_1+\|\Phi^T_2 S_j\|_1 \quad \text{subject to }
\,\, f_j = W_j + S_j.
\end{equation}
Notice that $\Phi^T_1 W_j$ and $\Phi^T_2 S_j$ are the wavelet and shearlet coefficients of the signals
$W_j$ and $S_j$, respectively. Notice that our objective is {\em not} on searching for the
sparsest expansion in a wavelet-shearlet dictionary, but on separation. Thus we can avoid an extensive,
presumably numerically instable search by minimizing specific coefficients, namely the {\em analysis} in contrast to the
{\em synthesis} coefficients, for each possible separation $f_j = W_j + S_j$.

We wish to further remark, that here the $\ell_1$ norm is placed
on the {\em analysis} rather than the {\em synthesis} coefficients to avoid numerical instabilities due to the
redundancy of the shearlet frame.

\subsection{Theoretical Result}

The theoretical result of the precision of separation of $f_j$ via (\ref{eq:sep}) proved in \cite{DK10} and \cite{DK09}
can now be stated in the following way:

\begin{theorem}[\cite{DK10} and \cite{DK09}]
\label{theo:geo_sep}
Let $\hat W_j$ and $\hat S_j$ be solutions to the optimization problem (\ref{eq:sep}) for each scale $j$.
Then we have
\[
\frac{\|\cP_j - \hat W_j\|_2+\|\cC_j - \hat S_j\|_2}{\|\cP_j\|_2+\|\cC_j\|_2} \rightarrow 0, \quad j \rightarrow \infty.
\]
\end{theorem}

This result shows that the components $\cP_j$ and $\cC_j$ are recovered with asymptotically arbitrarily high
precision at very fine scales. The energy in the pointlike component is completely captured by the wavelet coefficients,
and the curvelike component is completely contained in the shearlet coefficients. Thus, the theory evidences that the
Geometric Separation Problem can be satisfactorily solved by using a combined dictionary of wavelets and shearlets
and an appropriate $\ell_1$ minimization problem.

\subsection{Extensions}

Our numerical scheme for image separation, which we will present in detail in the next section, will use a shift invariant wavelet
tight frame and a compactly supported shearlet frame as opposed to orthonormal Meyer wavelets and band-limited shearlets required
for Theorem \ref{theo:geo_sep}. Hence this deserves some comments. Firstly, Theorem  \ref{theo:geo_sep} is based on an abstract
separation estimate which holds for any pair of frames, provided certain relative sparsity and cluster coherence
conditions with respect to the components of the data to be separated are satisfied (cf. \cite{DK10}). Secondly, using
the recently introduced concept of {\em sparsity equivalence} (see \cite{DK09,Kut10}), results requiring sparsity and
coherence conditions can be transferred from one system (set of systems) to another by `merely' considering particular
decay conditions of the cross-Grammian matrix (matrices). We strongly believe that this framework allows a similar
result as Theorem \ref{theo:geo_sep} for the pair of a shift invariant wavelet tight frame and a compactly supported shearlet frame.
Since the focus of this paper is however on the introduction of the numerical scheme, such a highly technical, theoretical analysis is
beyond the scope of this paper and will be treated in a subsequent work.

\section{Our Algorithmic Approach to the Geometric Separation Problem}
\label{sec:code}

In this section, we present our algorithmic approach to the Geometric Separation Problem of separating
point- from curvelike features by using a combined dictionary of wavelets and shearlets. The ingredients
of the algorithm will be detailed below.

\subsection{General Scheme}

In practice, the observed signal $f$ is often contaminated by noise which requires an adaption of the
optimization problem (\ref{eq:sep}). As proposed in numerous publications, one typically considers a modified
optimization problem -- so-called Basis Pursuit Denoising (BPDN) -- which can be obtained by relaxing
the constraint in (\ref{eq:sep}) in order to deal with noisy observed signals (see \cite{Ela10}). For
each scale $j$, the optimization problem then takes the form:
\begin{equation}\label{eq:BPDN}
(\hat W_j,\hat S_j) = \text{argmin}_{W_j,S_j}\|\Phi^T_1 W_j\|_1+\|\Phi^T_2 S_j\|_1+\lambda
\|f_j-W_j - S_j\|_2^2.
\end{equation}

In this new form, the additional content in the image -- the noise -- characterized by the property that it
can not be represented sparsely by either one of the two representation systems will be allocated
to the residual $f_j-W_j-S_j$. Hence, performing this minimization, we not only separate point- and curvelike
objects, which were modeled by $\cP_j$ and $\cC_j$ in Subsection \ref{subsec:model}, but also succeed in
removing an additive noise component as a by-product. Of course, solving the optimization problem (\ref{eq:BPDN}) for
all relevant scales $j$ is computationally expensive.

\subsection{Preprocessing}\label{subsec:pre}

To avoid high complexity, we observe that the frequency distribution of point- and curvelike components is
highly concentrated on high frequencies. Hence it would be essentially sufficient for achieving accurate
separation to solve (\ref{eq:BPDN}) for only sufficiently large scales $j$, as also evidenced by Theorem
\ref{theo:geo_sep}. This idea leads to a simplification of the problem (\ref{eq:BPDN}) by modifying the
observed signal $f$ as follows: We first consider bandpass filters $F_0, \dots, F_L$, where $(F_j)_{j=0,\dots,L}$
is the family of bandpass filters defined in Subsection \ref{subsec:filtering} up to scale $L$, and $F_0$ is a
lowpass filter. Thus the observed signal $f$ satisfies
\[
f = \sum_{j=0}^{L} F_j \ast f_j \quad \text{with} \,\, f_j = F_j \ast f.
\]
For each scale $j$, we now carefully choose a non-uniform weight $w_j>0$ satisfying, in particular, $w_j < w_{j'}$,
if $j < j'$. These weights are then utilized for a weighted reconstruction of $f$ resulting in a newly
constructed signal $\tilde f$ by computing
\begin{equation} \label{eq:reweighting}
\tilde f = \sum_{j=0}^{L} w_j \cdot (F_j \ast f_j).
\end{equation}
In this way, the two morphological components, namely points and curves, can be enhanced by suppressing the
low frequencies.

While emphasizing that certainly other weights can be applied in our scheme, for the numerical tests presented
in this paper we chose $L=3$, i.e., 4 subbands, and weights $w_0 = 0, w_1 = 0.1, w_2 = 0.7$, and $w_3 = 0.7$.
This coincides with the intuition that a strong weight should be assigned to a band on which the power spectrum of
the underlying morphological contents is highly concentrated. We do not claim that this is necessarily the optimal
choice, and a comprehensive mathematical optimality analysis is beyond our reach at this moment. To our mind, the
numerical results though justify this choice.

\subsection{Solver for the $\ell_1$ Minimization Problem}
\label{subsec:l1}

Using the reweighted reconstruction of $f$  from \eqref{eq:reweighting} coined $\tilde{f}$, we now consider
the following new minimization problem:
\begin{equation}\label{eq:MCA}
(\hat W,\hat S) = \text{argmin}_{W,S}\|\Phi^T_1 W\|_1+\|\Phi^T_2 S\|_1+\lambda \|\tilde f-W - S\|_2^2.
\end{equation}
Note that the frequency distribution of $\tilde f$ is highly concentrated on the high frequencies -- in
other words, scaling subbands of large scales $j$ --, and Theorem \ref{theo:geo_sep} justifies our
expectation of a very precise separation using $\tilde f$ instead of $f$. Even more advantageous, the reduced problem
(\ref{eq:MCA}) no longer involves different scales $j$, and hence can be efficiently solved by various
fast numerical schemes. In our separation scheme, we use the same optimization method as the one used in MCALab
to solve \eqref{eq:MCA} now applied to a combined wavelet-shearlet dictionary. We refer to \cite{FSED09} for a detailed
description of an algorithmic approach to solve \eqref{eq:MCA}; see also \cite{Ela10}.
In the following subsections, we discuss the particular form of the matrices $\Phi^T_1$ and $\Phi^T_2$ which encode the
wavelet and shearlet transform in the minimization problem (\ref{eq:MCA}) we aim to solve.

\subsection{Wavelet Transform}
Let us start with the wavelet transform.  The undecimated digital wavelet transform is certainly the most fitting version of the
wavelet transform for the filtering of data, and hence this is what we utilize also here. This transform is obtained by skipping
the subsampling, thereby yielding an overcomplete transform, which in addition is shift-invariant. The redundancy factor of this
transform is $3J+1$, where $J$ is the number of decomposition levels. We refrain from further details and merely refer the reader
to \cite{Mal98}.

\subsection{Shearlet Transform}
\label{subsec:shearlet}

For the shearlet transform, we employ the digital shearlet transform implemented by 2D convolution with discretized compactly
supported shearlets, which was introduced in \cite{Lim11}, see also \cite{KLZ11}. In the earlier work \cite{Lim09b},
an faithful digitalization of the continuum domain shearlet transform using compactly supported shearlets generated by
separable functions has been developed. However, firstly, this algorithmic realization allows only a limited directional
selectivity due to separability and, secondly, compactly supported shearlets generated by separable functions do not form
a tight frame which causes an additional computational effort to approximate the inverse of the shearlet transform by iterative
methods. These problems have been resolved in \cite{Lim11,KLZ11} by using non-separable compactly supported generators, and
we now summarize this procedure.

In the sequel, we will discuss the implementation strategy for computing the shearlet coefficients $\langle f,\psi_{j,k,m}\rangle$
only for $\psi_{j,k,m} \in \Psi(\psi,1)$. The same procedure can be applied to shearlets in $\tilde \Psi(\tilde \psi,1)$ except for
switching the role of variables. Without loss of generality, let $j/2$ be an integer. For $J>0$ fixed, assume that
$
f(x) = \sum_{n \in \Z^2} f_J(n)2^L\phi(2^Lx-n),
$
where $\phi$ is a 2D separable scaling function of the form $\phi_1(x_1)\phi_1(x_2)$ satisfying
\begin{equation}\label{eq:scale}
\phi(x) = \sum_{n \in \Z^2} h(n)2\phi(2x-n).
\end{equation}
Let $\psi$ be a 2D separable wavelet defined by $\psi(x_1,x_2)=\psi_1(x_1)\phi_1(x_2)$, where $\psi_1$ is a 1D wavelet.
Further, assume that $\psi$ can be written as
\begin{equation}\label{eq:wavelet}
\psi(x) = \sum_{n \in \Z^2} w(n) 2\phi(2x-n),
\end{equation}
where $w(n)$ are 2D separable wavelet filter coefficients associated with scaling filter
coefficients $h(n)$. For each $j \ge 0$, define the (non-separable) shearlet generator $\psi^{\text{non}}_j$ by
\begin{equation}\label{eq:newshear}
\hat{\psi}^{\text{non}}_j (\xi) = P_{J-j/2}(\xi)\hat \psi(\xi),
\end{equation}
where $P_{\ell}(\xi) = P(2^{\ell+1}\xi_1,\xi_2)$ for $\ell \ge 0$ and the trigonometric polynomial $P$ is
a 2D fan filter (c.f. \cite{DV05}).
To implement
$
\langle f(\cdot),\psi^{\text{non}}_{j,k,m}(\cdot)\rangle = \langle f(\cdot),\psi^{\text{non}}_{j,0,m}(S_{2^{-j/2}k}\cdot)\rangle,
$
we make two observations: Firstly, the functions $\psi^{\text{non}}_{j,0,m}$ are wavelets generated by refinement equations
\eqref{eq:scale}, \eqref{eq:wavelet} and \eqref{eq:newshear}. Thus, for each $j$, there exists an associated 2D wavelet filter
$w_j$. Secondly, the shear operator $S_{2^{-j/2}k}$ can be faithfully discretized by the digital shear operator $S^d_{2^{-j/2}k}$
(see \cite{Lim09b,Lim11}, also \cite{KLZ11}). The {\em digital (non-separable) shearlet transform} is then, using the shearlet filters
$\psi^d_{j,k} = S^d_{2^{-j/2}k}(w_j)$, defined by
$$
SH(f_J)(m_1,m_2) = (f_J * \psi^d_{j,k})(2^{J-j}m_1,2^{J-j/2}m_2) \quad \text{for}\,\, f_J \in \ell^2(\Z^2).
$$
If downsampling by $A_{2^j}$ is omitted, a shift invariant shearlet transform $(f_J * \psi^d_{j,k})(m_1,m_2)$ is obtained, in
which case dual shearlet filters $\tilde{\psi}^{\text{d}}_{j,k}$ can be easily computed by deconvolution.
We then obtain the reconstruction formula
\[
f_J = \sum_{j,k} (f_J* \psi^{\text{d}}_{j,k}(- \, \cdot))*\tilde{\psi}^{\text{d}}_{j,k}.
\]

\section{Numerical Results}
\label{sec:numerics}

In this section, we present and discuss some numerical results of our proposed scheme for separating point-
and curvelike features. In each experiment we compare our scheme, which is freely available in the
ShearLab\footnote{ShearLab (Version 1.1) is available from \url{http://www.shearlab.org}.} toolbox, with
the separation algorithm MCALab\footnote{MCALab (Version 120) is available from \url{http://jstarck.free.fr/jstarck/Home.html}.}.
In contrast to our algorithm, MCALab uses wavelets and curvelets to separate point- and curvelike components,
and we refer to \cite{FSED09} for more details on the algorithm.

\subsection{Comparison by Visual Perception}

\begin{figure}[h]
\begin{center}
\includegraphics[width=1.2in]{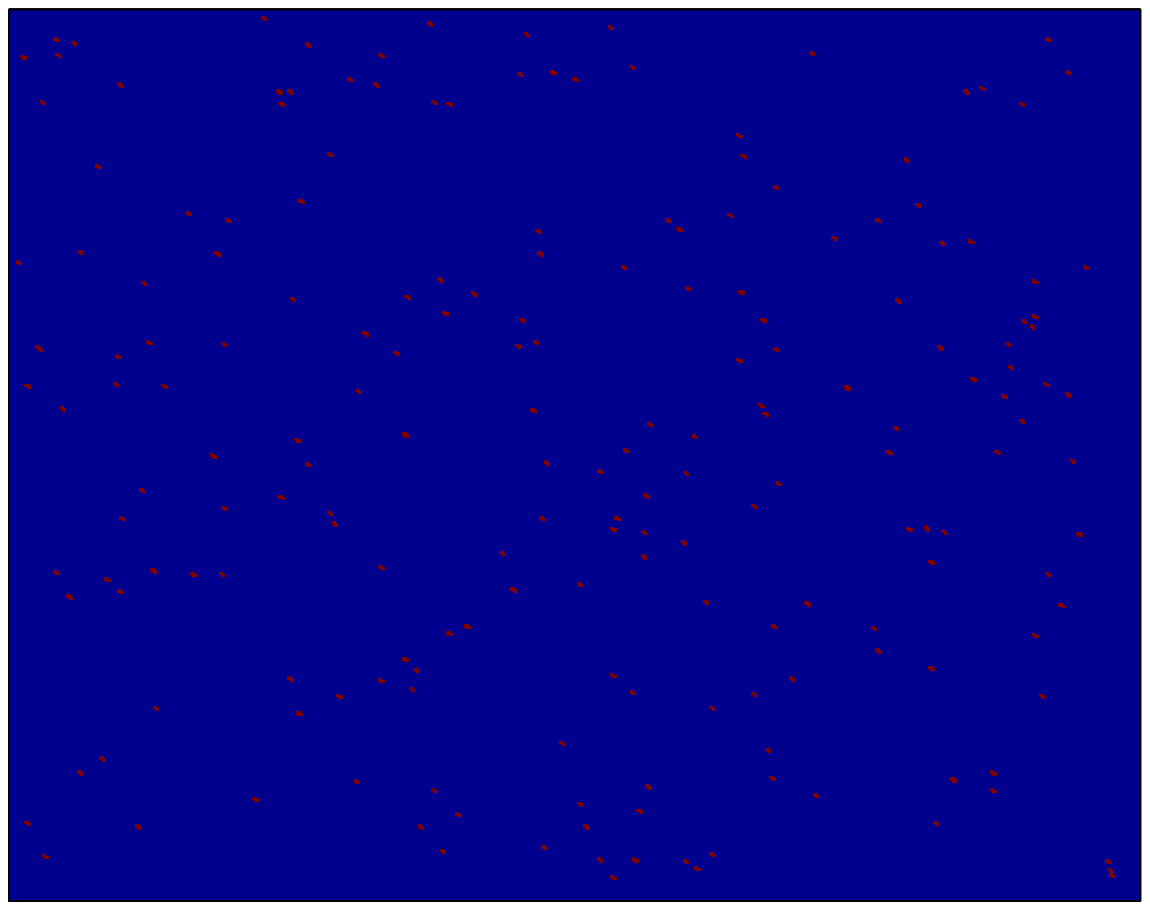}
\includegraphics[width=1.2in]{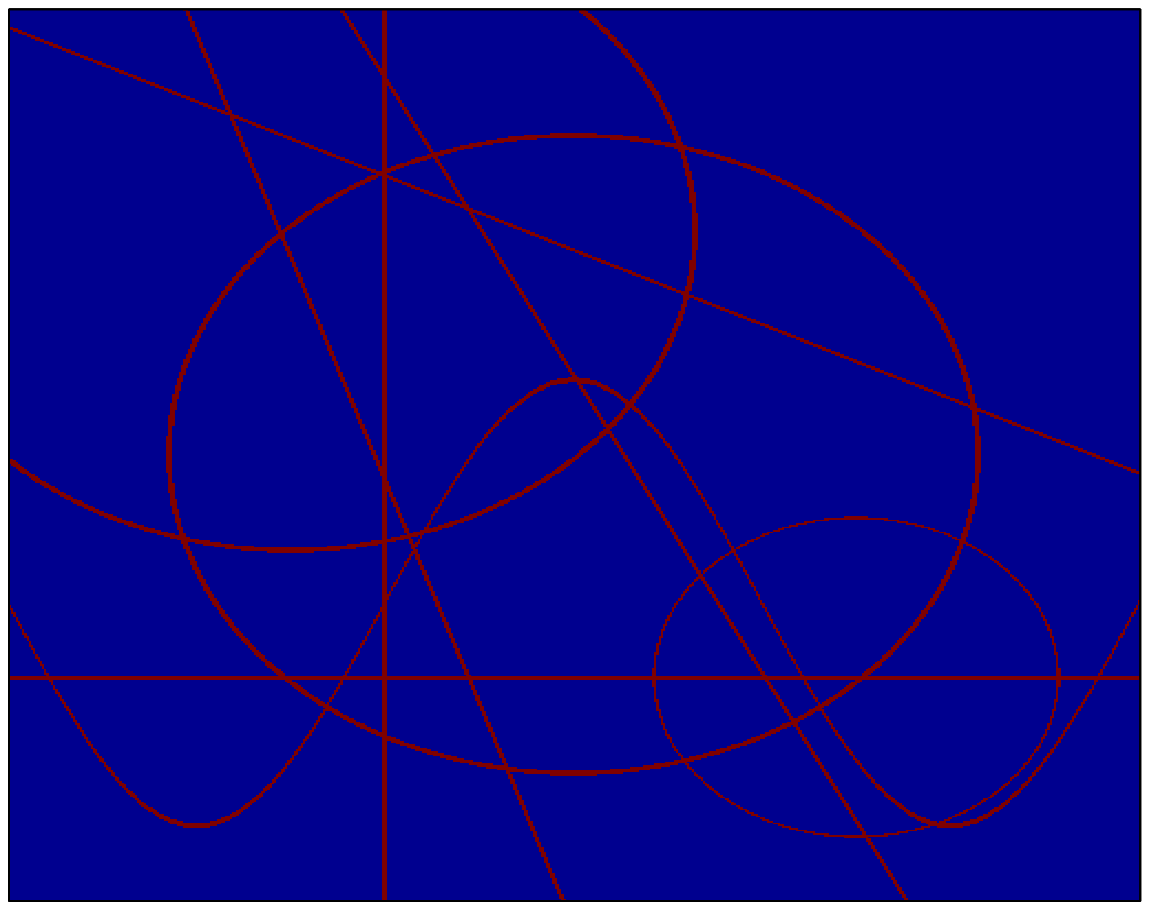}
\includegraphics[width=1.2in]{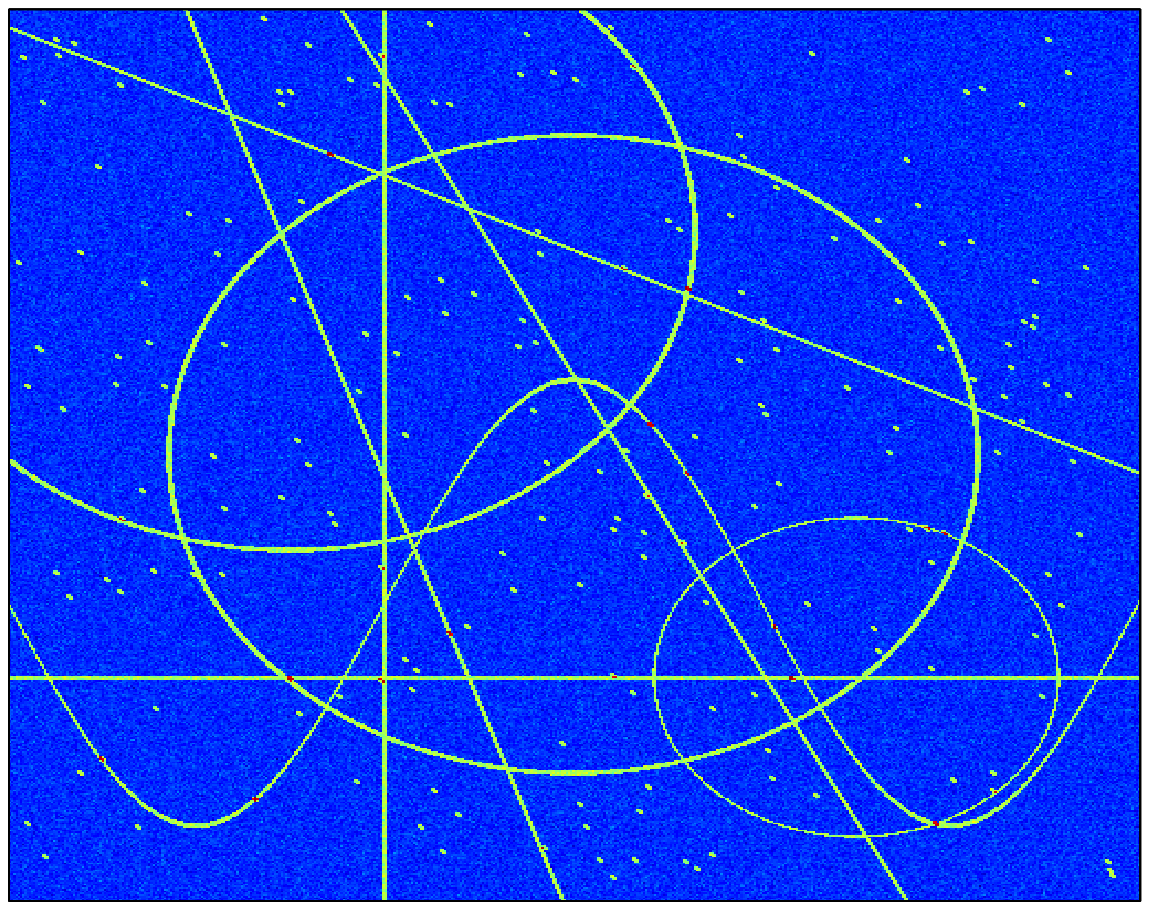}
\put(-45,-10){(c)}
\put(-140,-10){(b)}
\put(-225,-10){(a)}
\caption{(a) P: Image of points. (b) C: Image of curves. (c) Noisy image (512 $\times$ 512).}
\label{fig:orig1}
\end{center}
\end{figure}

Aiming first at comparison by visual perception, we choose an artificial image $I$ composed of two subimages
$P$ and $C$, where $P$ solely contains pointlike structures and $C$ different curvelike structures
(see Figures \ref{fig:orig1}(a) and (b)). We further add white Gaussian noise to $I = P+C$, shown in Figure
\ref{fig:orig1}(c).

\begin{figure}[h]
\begin{center}
\includegraphics[width=1.5in]{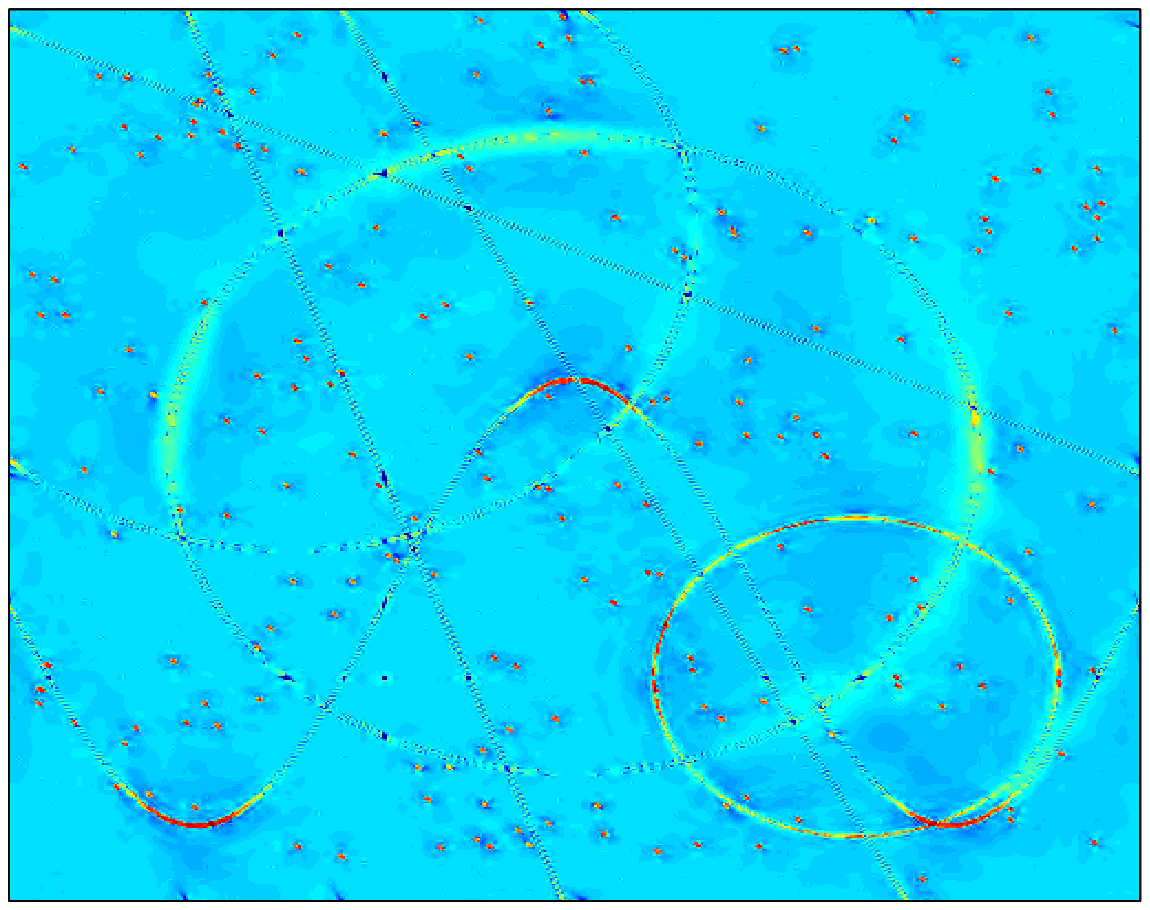}
\includegraphics[width=1.5in]{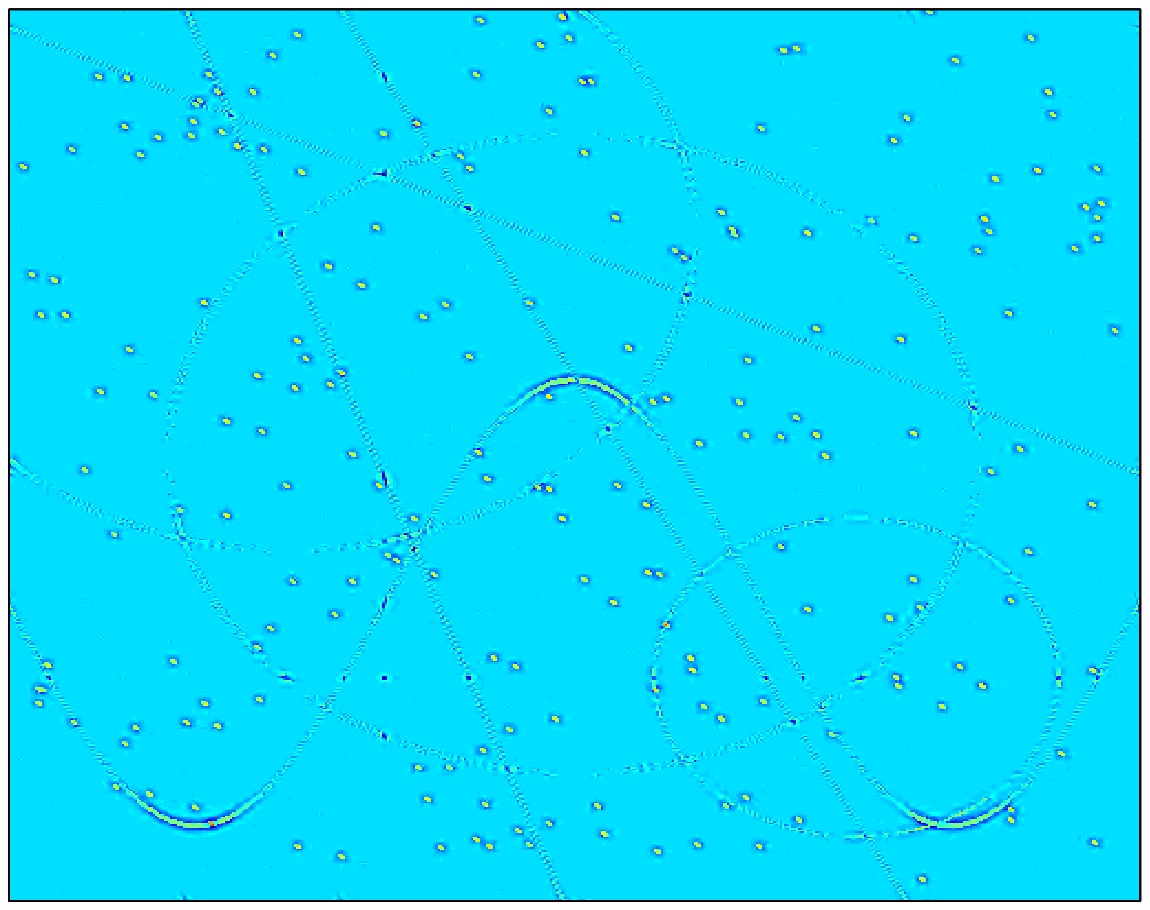}
\put(-170,0){(a)}
\put(-60,0){(b)}
\caption{ (a) Pointlike image component extracted by MCALab without preprocessing. (b) Pointlike image component extracted by  MCALab with preprocessing.}
\label{fig:band}
\end{center}
\end{figure}

Our scheme consists of two parts: Preprocessing of the image as described in Subsection \ref{subsec:pre}, followed by separation
using a combined wavelet-shearlet dictionary as described in Subsections \ref{subsec:l1}-\ref{subsec:shearlet}. First,
we will focus on the preprocessing step, and apply MCALab
with and without our preprocessing step -- due to limited space we just mention that our scheme is similarly
positively affected by preprocessing. Figures \ref{fig:band}(a) and (b) then visually indicate that preprocessing indeed
significantly improves the accuracy of separation.

\begin{figure}[h]
\begin{center}
\includegraphics[width=1.5in]{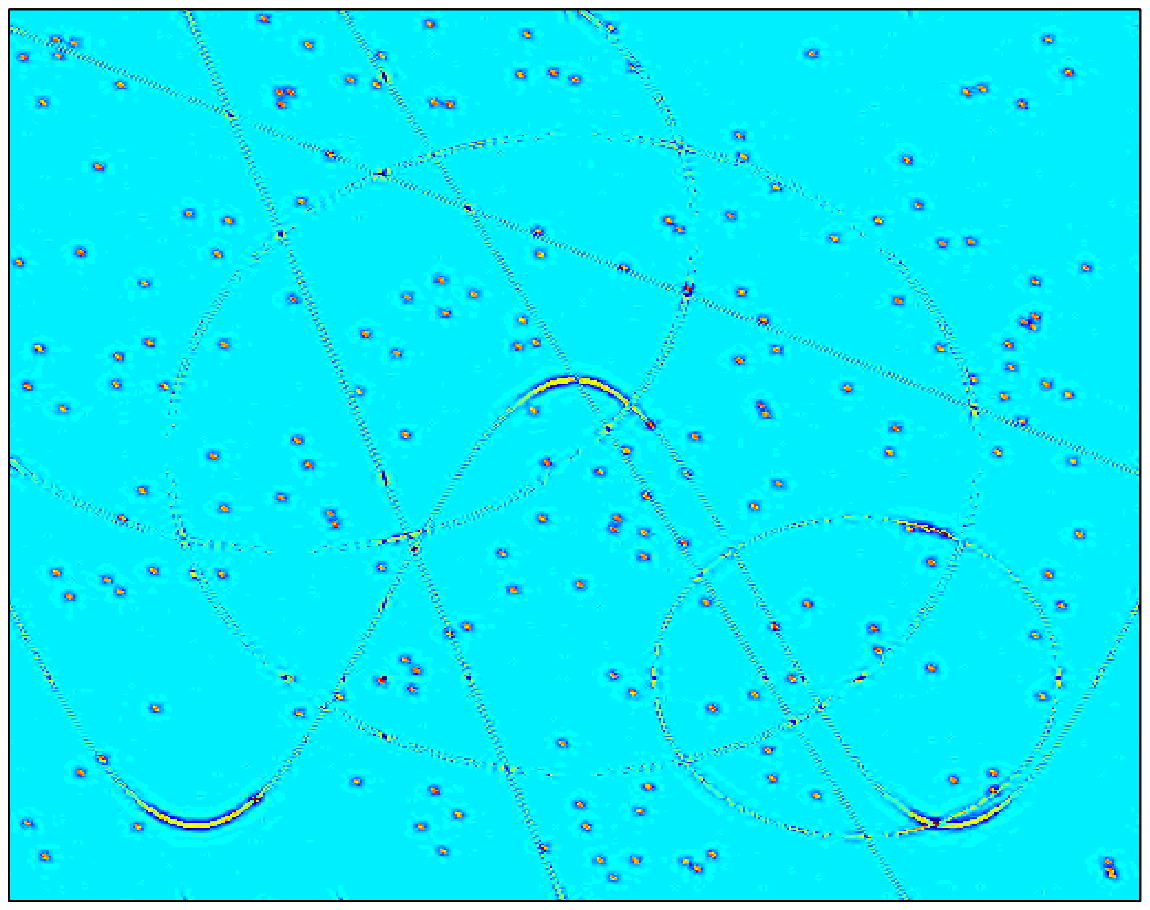}
\includegraphics[width=1.5in]{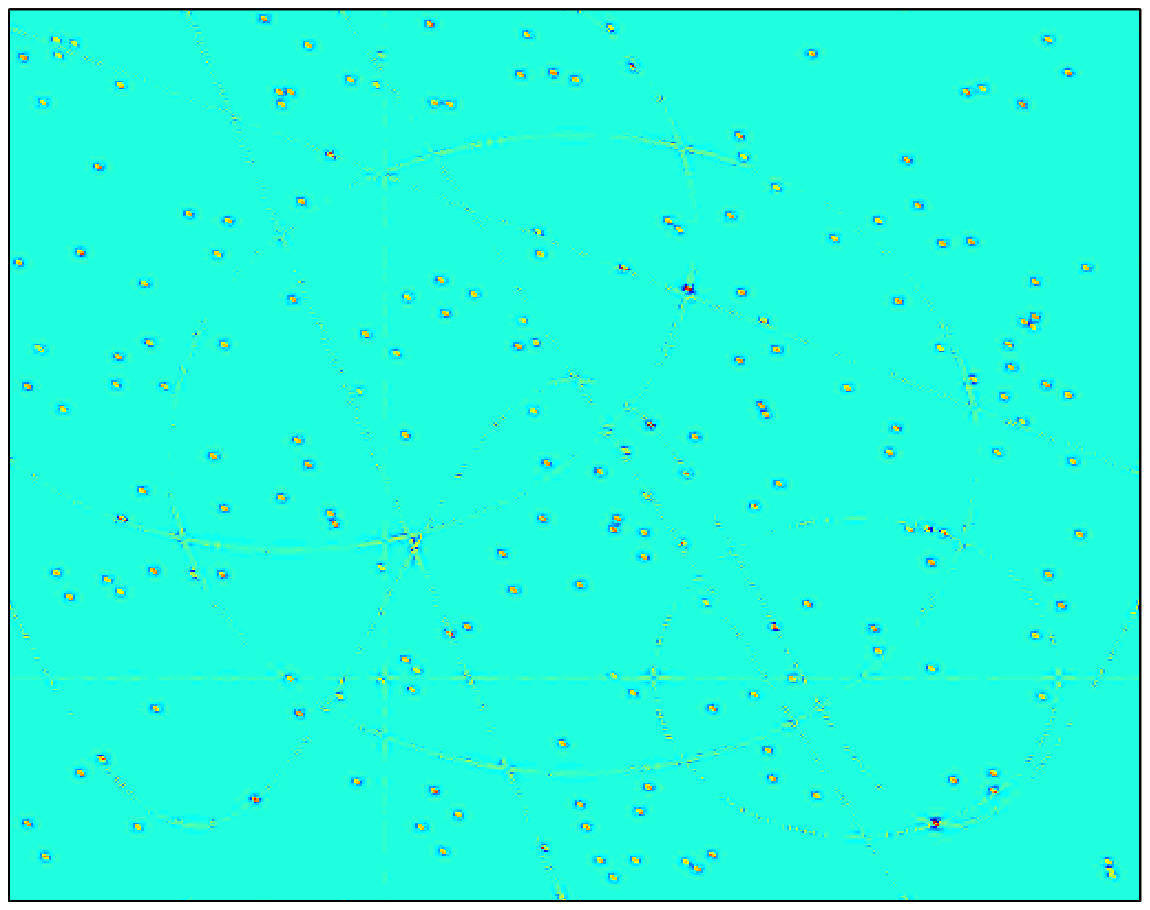}
\put(-170,0){(a)}
\put(-60,0){(b)}
\\
\includegraphics[width=1.5in]{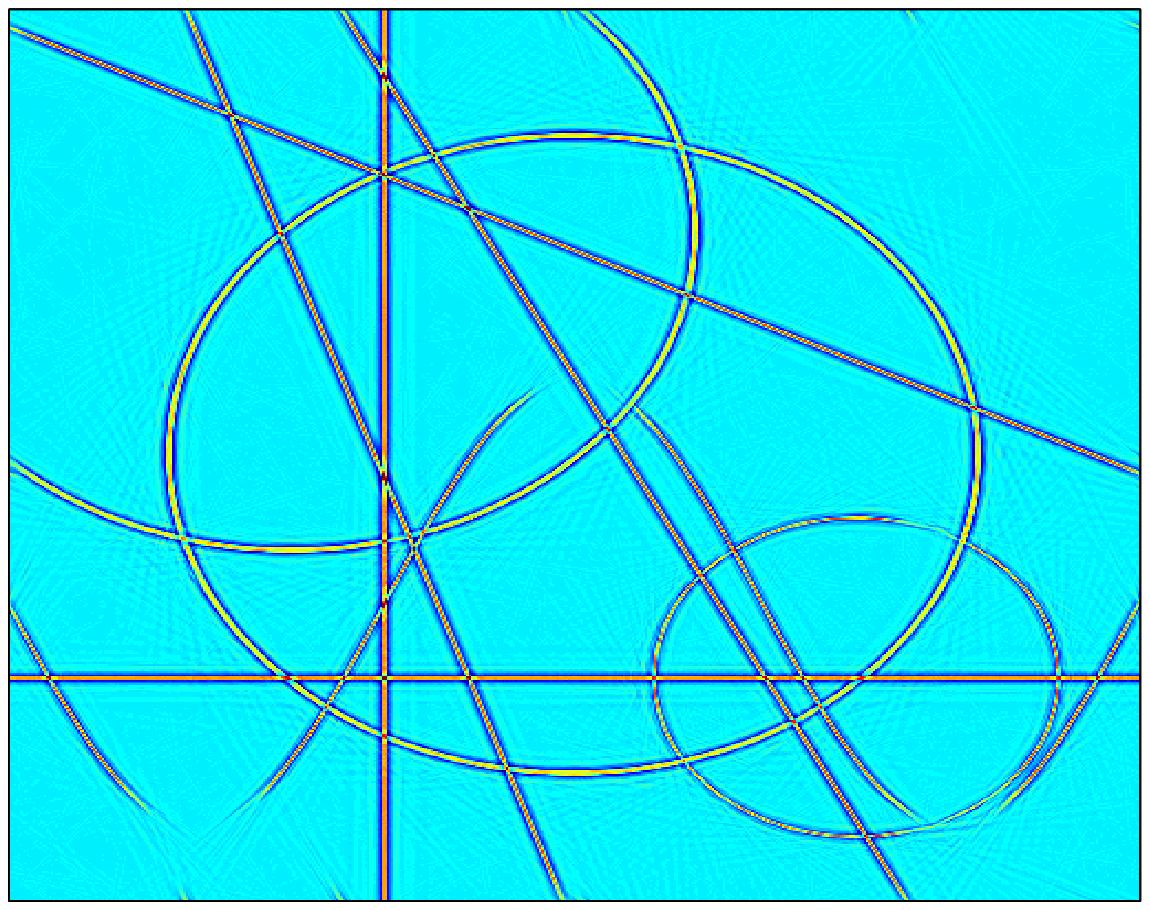}
\includegraphics[width=1.5in]{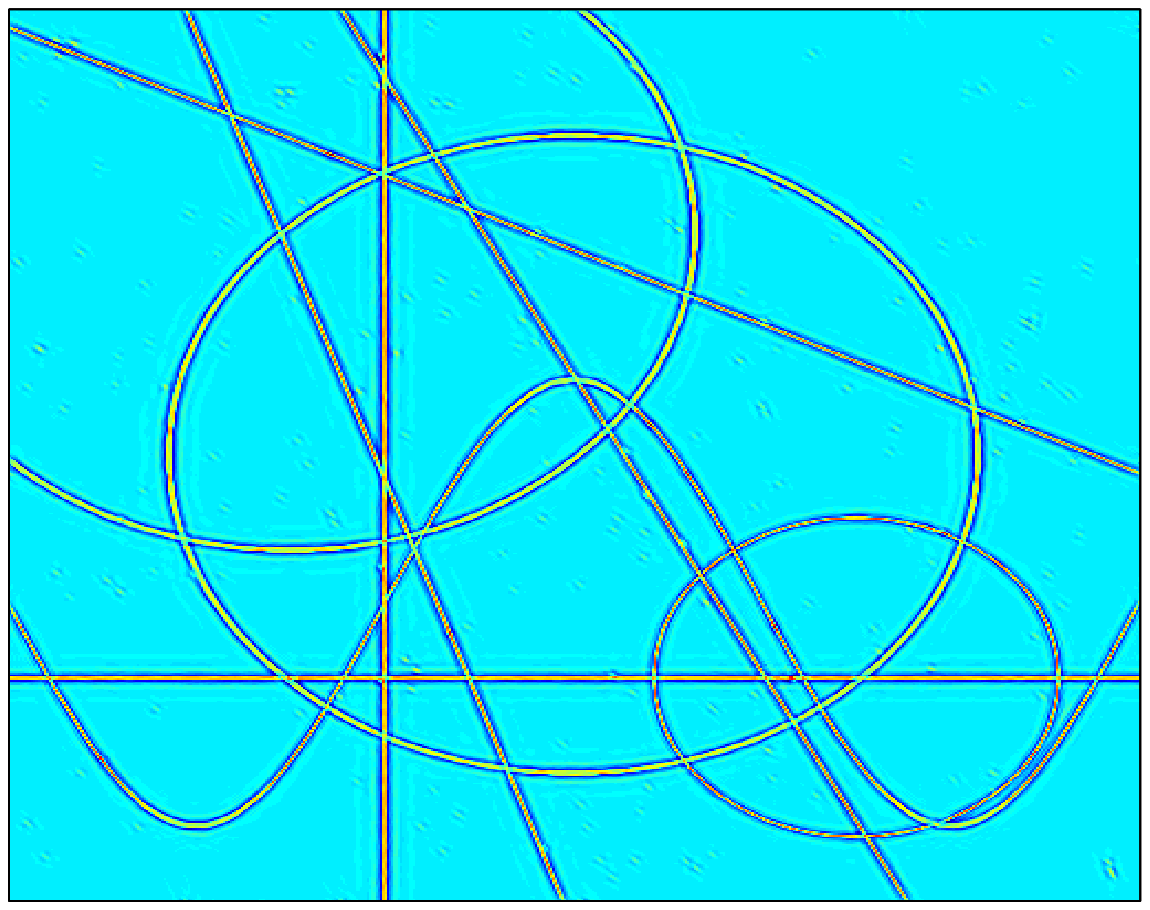}
\put(-170,0){(c)}
\put(-60,0){(d)}
\caption{ (a) MCALab: Pointlike component.  (b) Our scheme: Pointlike component. (c) MCALab: Curvelike component. (d) Our scheme: Curvelike component.}
\label{fig:comp}
\end{center}
\end{figure}

\begin{figure}[h]
\begin{center}
\includegraphics[width=1.5in]{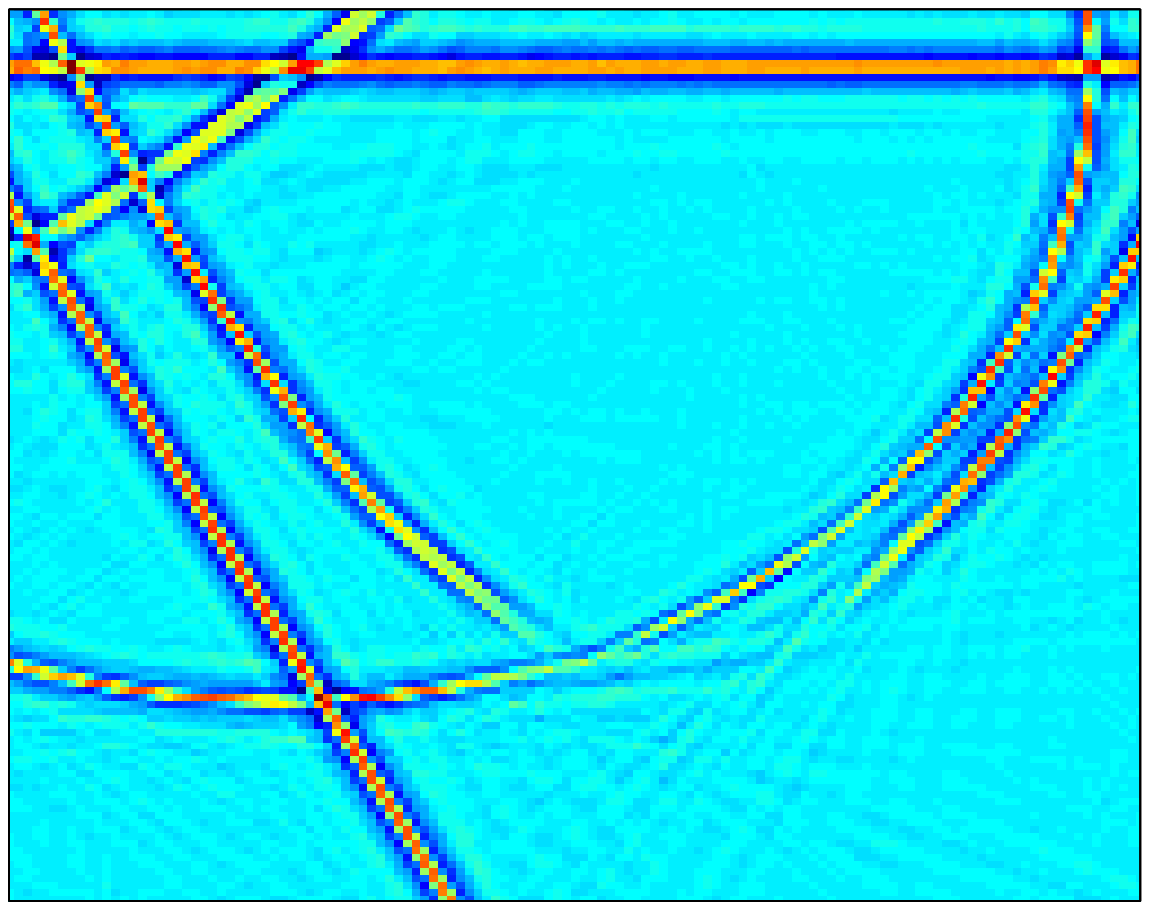}
\includegraphics[width=1.5in]{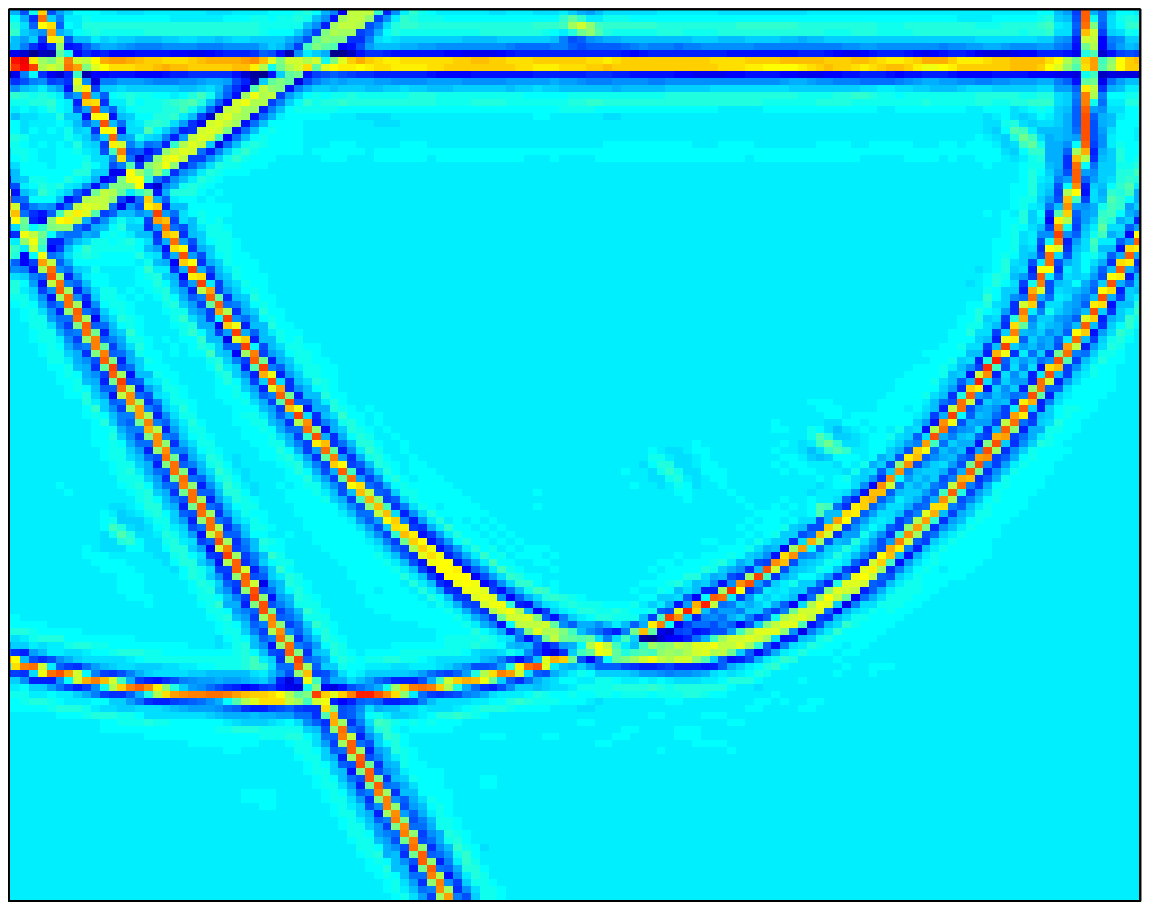}
\put(-170,0){(a)}
\put(-60,0){(b)}
\caption{ Zoomed images: (a) MCALab.  (b) Our scheme.}
\label{fig:zoom}
\end{center}
\end{figure}

To achieve a fair comparison, we now apply both schemes to the same {\em preprocessed} image as defined in \eqref{eq:reweighting}.
Figures \ref{fig:comp} (a)-(d) and Figures \ref{fig:zoom}(a)-(b) show the comparison results. In Figure \ref{fig:comp}(c), it can be
observed that the curvelet transform performs well for extracting lines due to excellent directional selectivity. However,
some part of the curve is missed and appears in the pointlike part, see also Figure \ref{fig:zoom}(a). This error becomes worse with growing curvature. In contrast
to this, compactly supported shearlets provide much better spatial localization than (band-limited) curvelets, which positively
affects the capturing of localized features of the curve as illustrated in Figure \ref{fig:comp}(d) and Figure \ref{fig:zoom}(b).
Hence, with respect to this visual comparison, our scheme outperforms MCALab, in particular, when curves with large
curvature are present.
\vspace{-20pt}
\vspace{10pt}

\subsection{Comparison by Quantitative Measures}

We now put the comparison on more solid ground by introducing two quantitative measures for analyzing how accurate
our scheme as compared to MCALab extracts points and curves.

\begin{figure}[h]
\begin{center}
\includegraphics[width=1.5in]{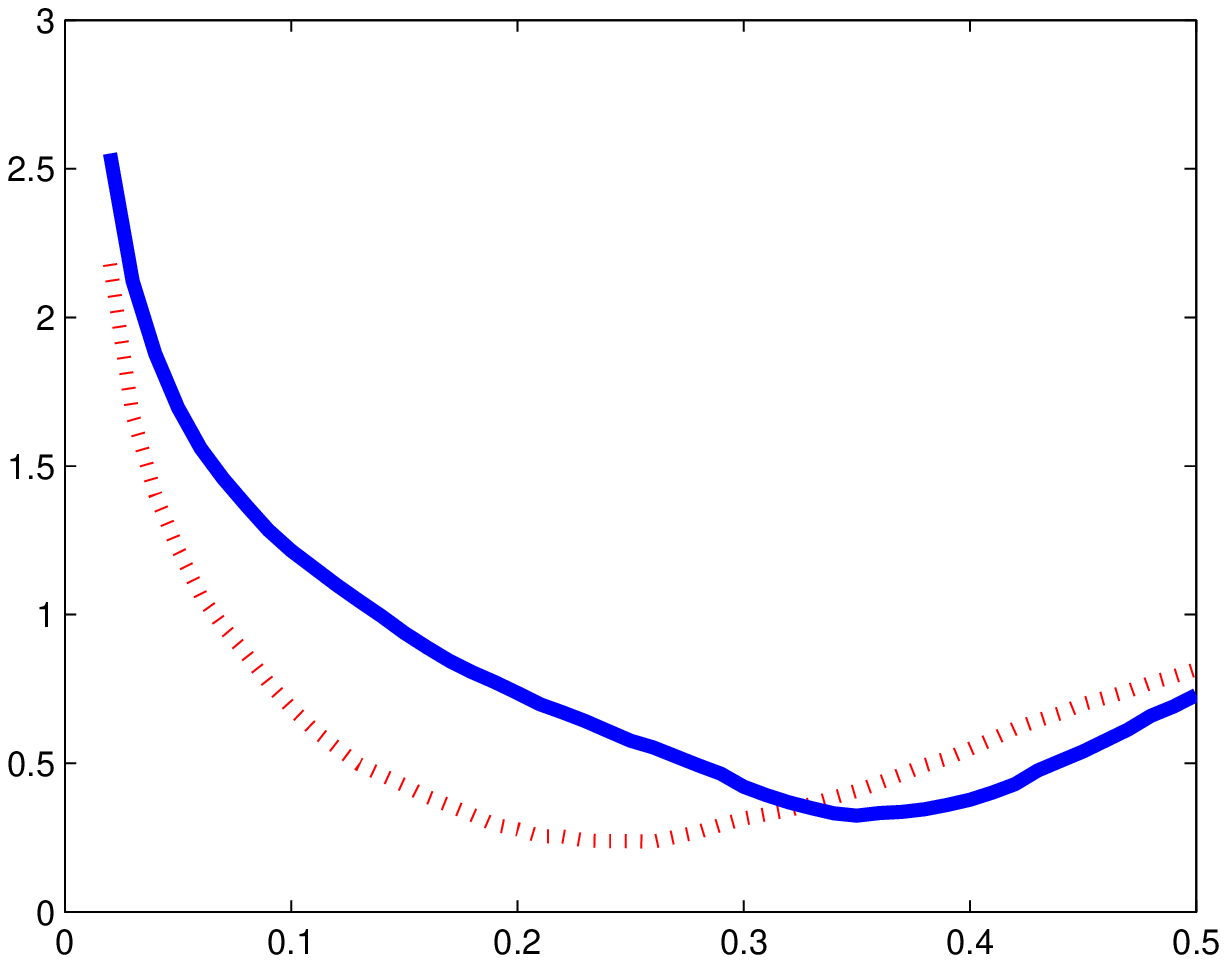}
\includegraphics[width=1.5in]{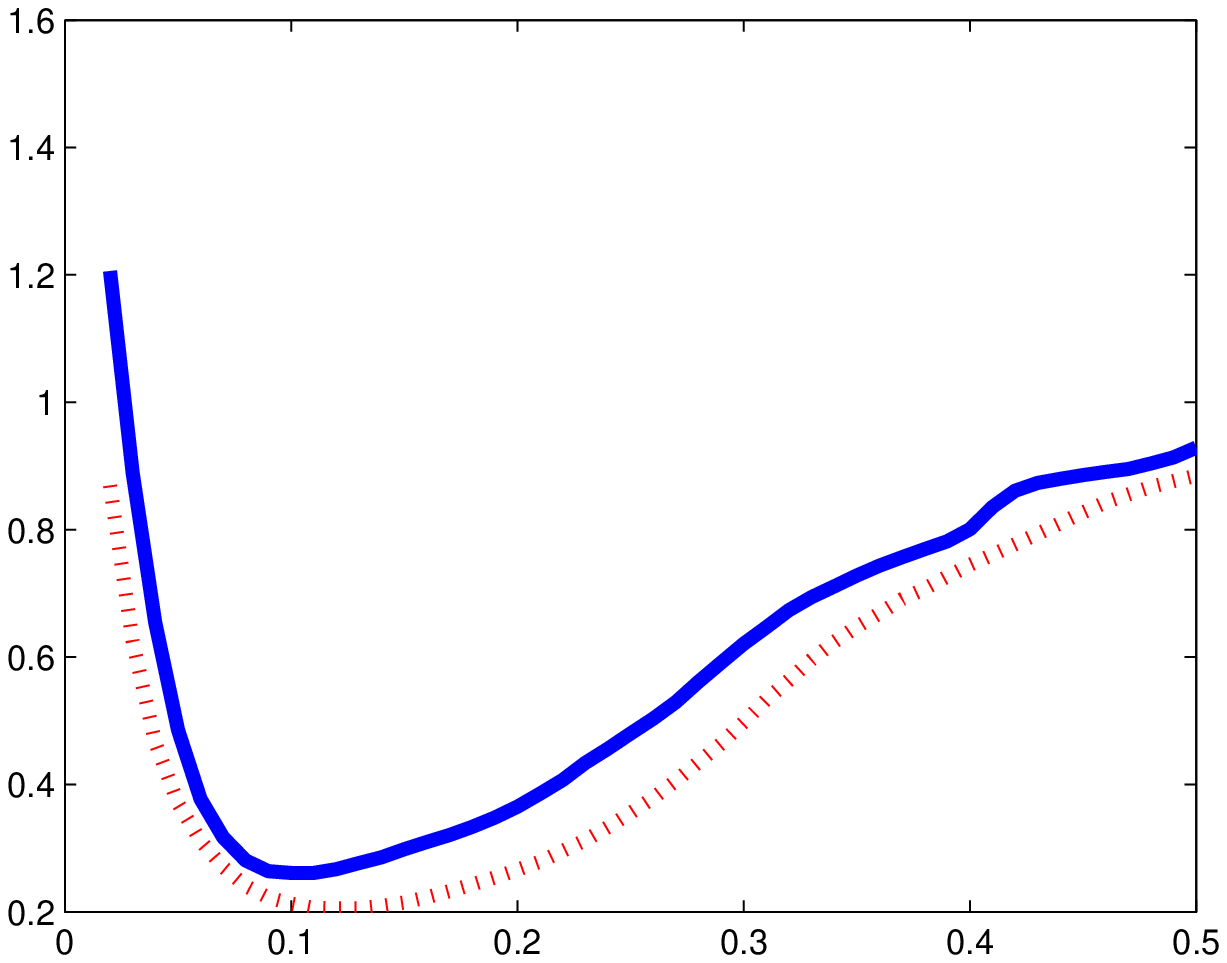}
\put(-170,0){(a)}
\put(-60,0){(b)}
\caption{(a) Graph of quantitative measure $T \mapsto M_p(\hat P)(T)$: Our scheme (dashed curve) and MCALab (solid curve).
(b) Graph of  quantitative measure $T \mapsto M_p(\hat C)(T)$: Our scheme (dashed curve) and MCALab (solid curve).}
\label{fig:error}
\end{center}
\end{figure}

To define our first quantitative measure, let $P$ and $C$ be (binary) images containing
points and curves, respectively, with image domain $\Omega \subset \Z^2$, and let $\hat P$ and $\hat C$ be the
separated images from $I = P+C+\text{noise}$ by the separation scheme to be analyzed. Letting $T \ge 0$, $\text{B}_T$ be
defined by
$$
\text{B}_T(I) = \chi_{\{n \in \Omega : I(n) \ge T\}}, \quad \text{for a 2D image}\,\, I,
$$
and $g$ be a 2D discrete Gaussian filter, we introduce the test measures
$$
M_p(\hat P)(T) = \frac{\|g*P-g*(\text{B}_{T\cdot\max(\hat P)}(\hat P))\|_2}{\|g*P\|_2}
$$
and
$$
M_c(\hat C)(T) = \frac{\|g*C-g*(\text{B}_{T\cdot\max(\hat C)}(\hat C))\|_2}{\|g*C\|_2}.
$$
Using $P$ and $C$ as given by Figure \ref{fig:orig1}, the graphs of the error functions $M_p(\hat P)$ and $M_c(\hat C)$
for our scheme and MCALab are plotted in Figure \ref{fig:error} depending on the threshold parameter $0<T<1$.
These figures imply that our scheme outperforms MCALab with respect to this quantitative measure.

As our second quantitative measure, we will use the running time of each scheme. With respect to this comparison measure, MCALab runs
182.19 sec (30 iterations) to produce the test results while our scheme takes 135.37 sec (15 iterations). The running time
was computed by taking the average over 10 runs. Again, with respect to this measure our scheme performs superior.

\section{Application in Neurobiology}

\begin{figure}[h]
\begin{center}
\includegraphics[width=1.2in]{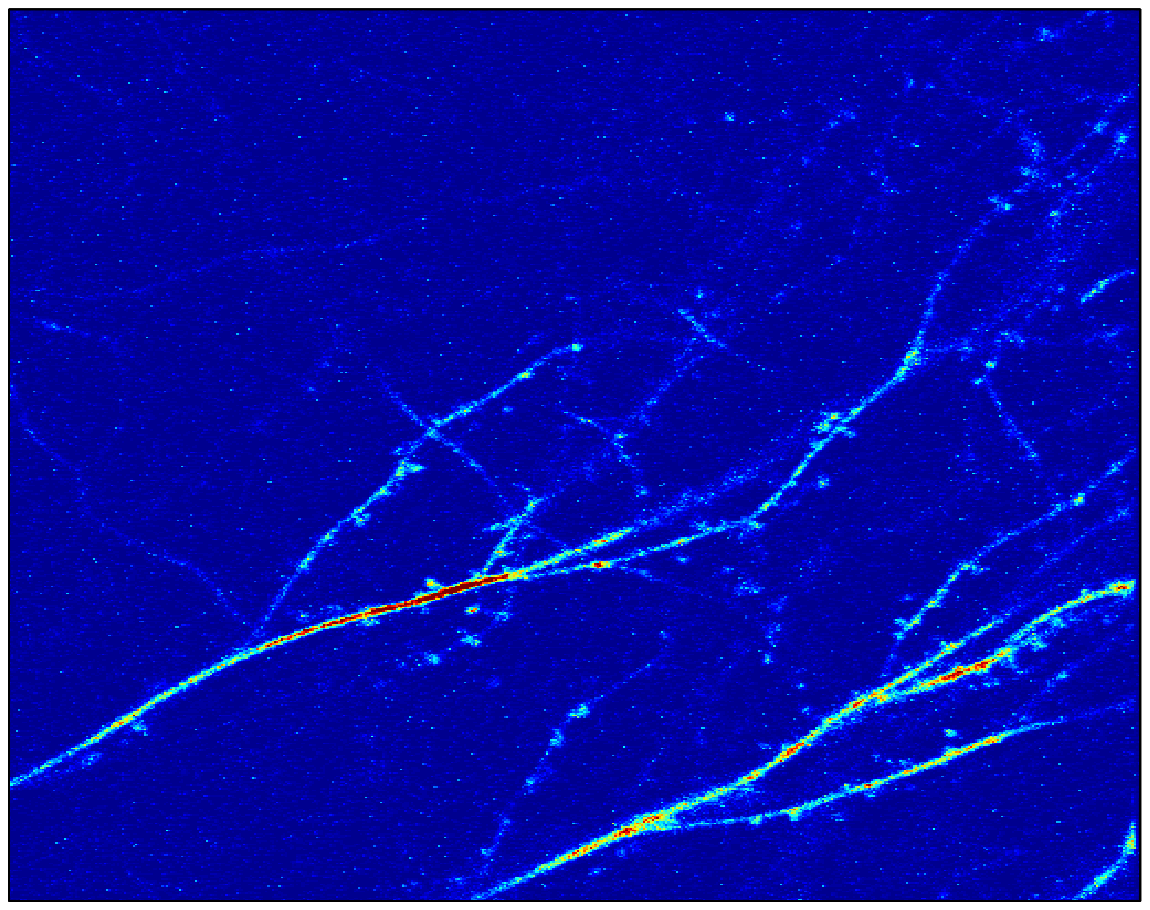}
\includegraphics[width=1.2in]{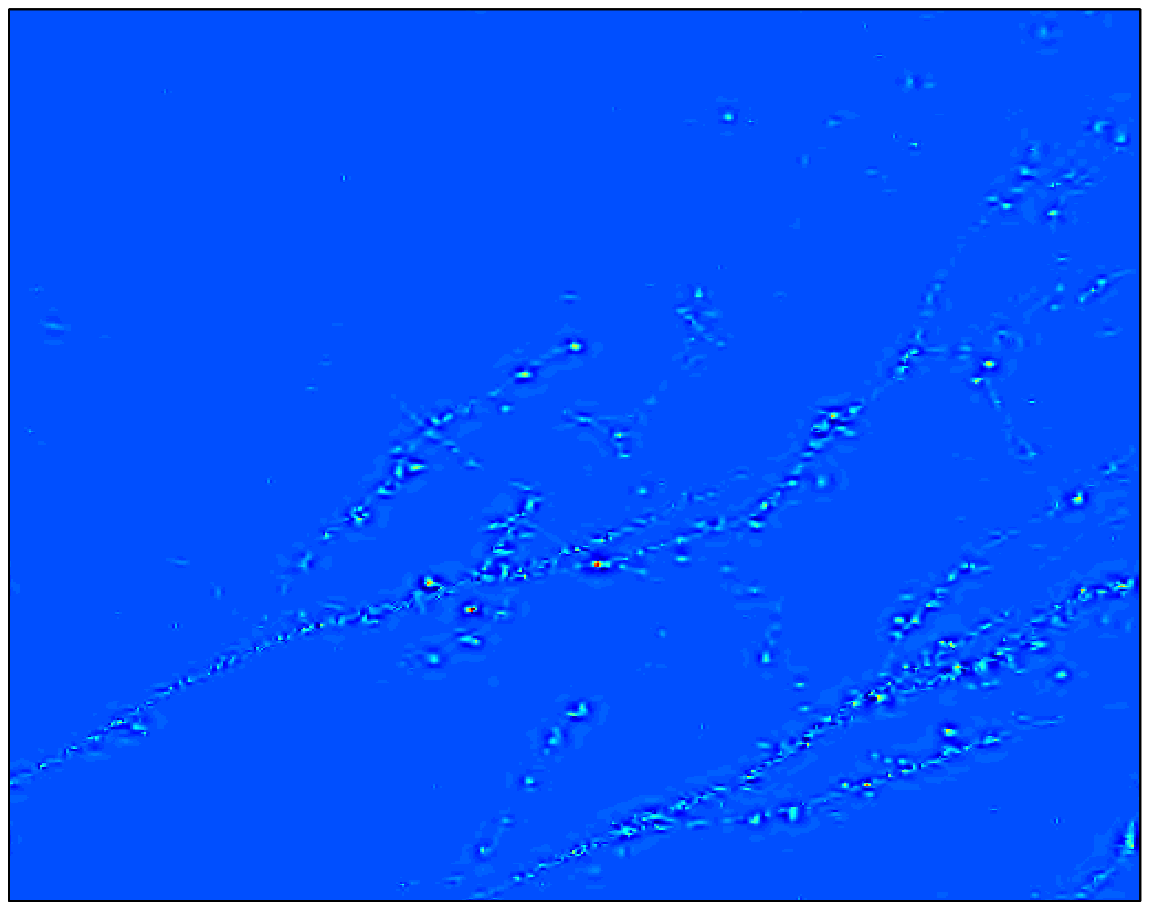}
\includegraphics[width=1.2in]{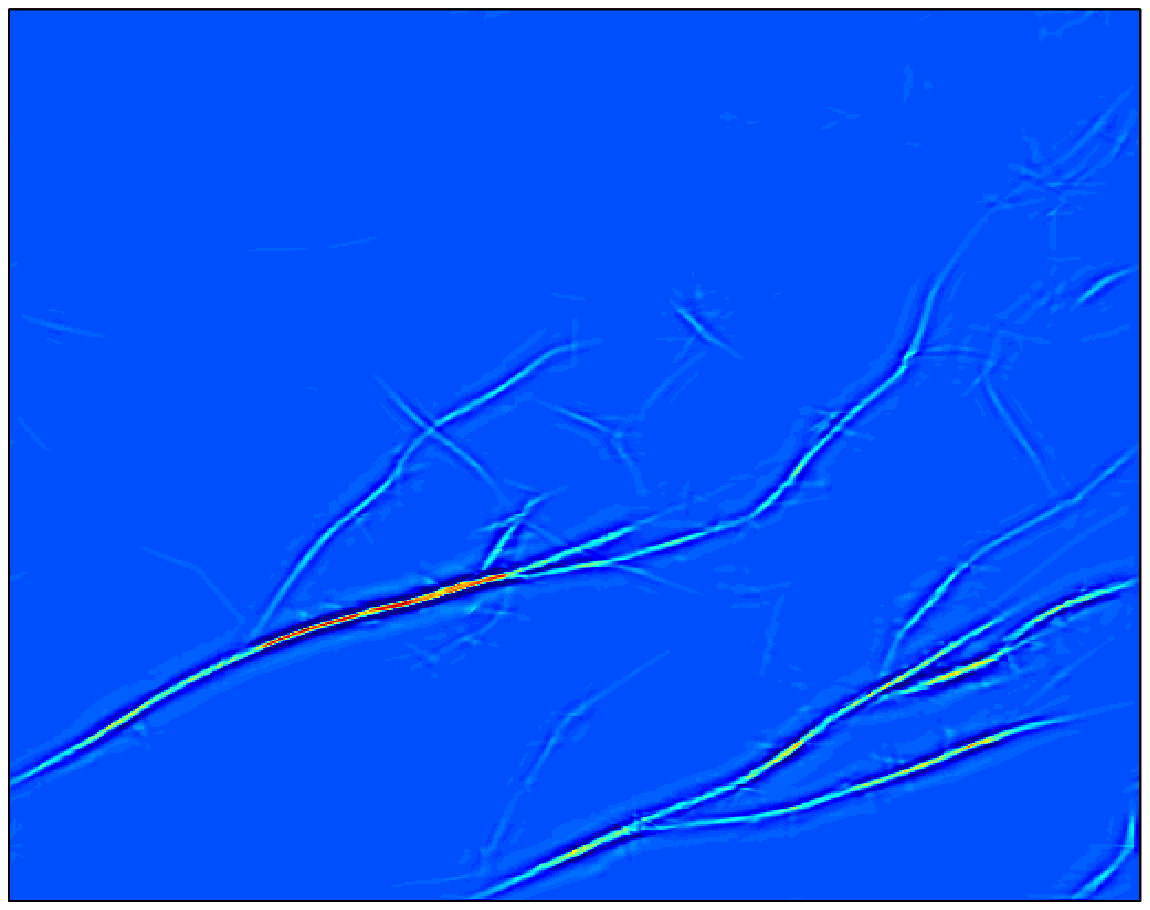}
\put(-45,-10){(c)}
\put(-140,-10){(b)}
\put(-225,-10){(a)}
\caption{(a) Image of neuron. (b) Extracted spines. (c) Extracted dendrites.}
\end{center}
\label{fig:neuron}
\end{figure}

To test the performance of our scheme on real-world images, we apply it to an image of a neuron generated by fluorescence microscopy
(Figure \ref{fig:neuron}(a)),  which is composed of `spines' (pointlike features) and `dendrites' (curvelike features).
Figures \ref{fig:neuron}(b) and (c) show the extracted images containing spines and dendrites, respectively.

\section{Conclusion}
\label{sec:conclusion}
In this paper, we introduced a novel methodology for separating images into point- and curvelike
features based on the new paradigm of sparse approximation and using $\ell_1$ minimization. In contrast to other approaches, our
algorithm utilizes a combined dictionary consisting of wavelets and shearlets, implemented as shift-invariant
transforms, and is based on a mathematical theory. The excellent localization property of compactly supported
shearlets allows shearlets to capture the curvelinear part very accurately and efficiently, even if the curvature is
relatively large. Numerical results show that our scheme
extracts point- and curvelike features more precise and uses less computing time than the
state-of-the-art algorithm MCALab, which is based on wavelets and curvelets.

\section*{Acknowledgement}
The first author would like to thank David Donoho and Michael Elad for inspiring discussions on related topics.
We are also grateful to the research group by Roland Brandt for supplying the test image in Figure \ref{fig:neuron}.
The authors acknowledge support from DFG Grant SPP-1324, KU 1446/13. The first author also acknowledges partial
support from DFG Grant KU 1446/14. We are particularly grateful to the two anonymous referees
for their many valuable comments and suggestions which significantly improved this paper.

\end{document}